\documentstyle [12pt]{article}

\def\qmod#1#2{{\hbox{}^{\displaystyle{#1}}}\!\big/\!\hbox{}_{\displaystyle{#2}}}

\newfam\msbfam
\font\twelmsb=msbm10 at 12pt
\font\tenmsb=msbm10
\font\sevenmsb=msbm10 at 7pt
\font\fivemsb=msbm10 at 5pt

\textfont\msbfam=\twelmsb
\scriptfont\msbfam=\sevenmsb
\scriptscriptfont\msbfam=\fivemsb

\def\Bbb{\fam\msbfam\tenmsb}

\def\h{{\germ h}}

\def\p{{\germ p}}

\def\C{{\Bbb C}}

\def\P{{\Bbb P}}

\def\R{{\Bbb R}}
\def\Z{{\Bbb Z}}

\def\qed {\hfill\vrule height6pt width6pt depth0pt \bigskip}

\def\map{\longrightarrow}
\def\textmap#1{\mathop{\vbox{\ialign{
				##\crcr
    ${\scriptstyle\hfil\;\;#1\;\;\hfil}$\crcr
    \noalign{\kern-1pt\nointerlineskip}
    \rightarrowfill\crcr}}\;}}

\def\textlmap#1{\mathop{\vbox{\ialign{
				##\crcr
    ${\scriptstyle\hfil\;\;#1\;\;\hfil}$\crcr
    \noalign{\kern-1pt\nointerlineskip}
    \leftarrowfill\crcr}}\;}}

\newfam\meuffam
\font\tenmeuf=eufm10
\font\sevenmeuf=eufm7
\font\fivemeuf=eufm5

\textfont\meuffam=\tenmeuf
\scriptfont\meuffam=\sevenmeuf
\scriptscriptfont\meuffam=\fivemeuf

\def\germ{\fam\meuffam\tenmeuf}

\begin{document}
\def\Pr{{\rm Pr}}
\def\tr{{\rm Tr}}
\def\End{{\rm End}}
\def\Pic{{\rm Pic}}
\def\NS{{\rm NS}}
\def\deg{{\rm deg}}
\def\det{{\rm det}}
\def\Hom{{\rm Hom}}
\def\Herm{{\rm Herm}}
\def\Vol{{\rm Vol}}
\def\pf{{\bf Proof: }}
\def\id{{\rm id}}
\def\ad{{\rm ad}}
\def\i{{\germ i}}
\def\im{{\rm im}}
\def\rk{{\rm rk}}
\def\Spin{{\rm Spin}}
\def\h{{\bf H}}
\def\U({{\rm U(}}
\def\SU{{\rm SU}}
\def\SO{{\rm SO}}
\def\dv{\bar\partial}
\def\dva{\bar\partial_A}
\def\da{\partial_A}
\def\p{\partial\bar\partial}
\def\pa{\partial_A\bar\partial_A}
\def\Dr{\hskip 4pt{\not}{D}}
\newtheorem{sz}{Satz}[section]
\newtheorem{th}[sz]{Theorem}
\newtheorem{pr}[sz]{Proposition}
\newtheorem{re}[sz]{Remark}
\newtheorem{co}[sz]{Corollary}
\newtheorem{dt}[sz]{Definition}
\newtheorem{lm}[sz]{Lemma}
\title{Seiberg-Witten Invariants and Rationality of Complex Surfaces}
\author{Christian Okonek$^*$\and Andrei Teleman\thanks{Partially supported
by: AGE-Algebraic Geometry in Europe, contract No ERBCHRXCT940557 (BBW
93.0187), and
by  SNF, nr. 21-36111.92} }
\date{}
\maketitle

\setcounter{section}{-1}
\section{Introduction}

Recently, Seiberg and Witten introduced new differential invariants for
4-manifolds, which are defined by counting solutions of the so called {\sl
monopole
equations}, a system of non-linear differential equations  of Yang-Mills-Higgs
type [18].

The new invariants are expected to be equivalent to the Donaldson  polynomial
invariants, and they have already found important applications [15].

The purpose  of this paper is:\\
- to explain the Seiberg-Witten invariants\\
- to show that --- on a K\"ahler surface --- the solutions of the monopole
equations can be interpreted as algebraic objects, namely effective divisors\\
- to give --- as an application --- a short selfcontained proof for the
fact that
rationality of complex surfaces is a ${\cal C}^{\infty}$-property.

\section{$\Spin^c$-structures and the monopole equation}

\begin{dt} {\rm [1], [11]} The group $\Spin^c(n):=\Spin(n)\times_{\Z_2}S^1$ is
called the complex spinor group.
\end{dt}

For the case $n=4$, there is a natural identification
$$\Spin^c(4)=\{(A,B)\in\U(2)\times\U(2)|\ \det A=\det B\}\ .$$

The following diagram summarizes some of the basic relations of
$\Spin^c(4)$ to other
groups:
$$\begin{array}{rrclcrl}
&&&&\U(2)&&\\
&&&^l\swarrow& &\searrow^i\\

\ \ \ S^1&\map&\ \ \Spin^c(4)\ &&\map&&\SO(4)\\
(\cdot)^2\downarrow\ \ &{\scriptstyle\det}\swarrow&\ \downarrow&&
\phantom{{\scriptstyle\lambda^+}}&&\ \ \ \ \ \downarrow
{\scriptstyle(\lambda^+,\lambda^-)}\\
\ \ \ S^1&&\ \U(2)\times U(2)&&\textmap{\ad}&&\SO(3)\times SO(3) \ \ \ \
\end{array}$$
Here \  $l:\U(2)\map\Spin^c(4)$ \ is the canonical lifting of the
homomorphism  \
$i\times\det:\U(2)\map\SO(4)\times S^1$ [11], and acts by the formula\linebreak
$\U(2)\ni a\longmapsto\left(\left(\matrix{\id&0\cr0&\det
a\cr}\right),a\right)\in\Spin^c(4)$.
$\lambda^{\pm}:\SO(4)\map\SO(3)$ are the maps induced by the two projections of
$\Spin(4)=\SU(2)^+\times\SU(2)^-$ onto the factors.

Let $X$ be a closed, oriented {\sl simply connected} 4-manifold,
$\Lambda^p$ the bundle of $p$-forms on $X$, and $A^p:=A^0(X,\Lambda^p)$ the
space of
sections in this bundle. Let $g$ be a Riemannian metric on $X$, denote by
$P$ the
associated principal
$\SO(4)$-bundle, and by $P^{\pm}$ the $\SO(3)$-bundles induced via the
morphisms
$\lambda^{\pm}$. The  real 3-vector bundles
$\Lambda^2_{\pm}:=P^{\pm}\times_{\SO(3)}\R^3$ can be identified with the
bundles of
(anti)self-dual 2-forms, hence there is an orthogonal splitting
$\Lambda^2=\Lambda^2_+\oplus\Lambda^2_-$.

\begin{lm}{\rm [10]} Given $c\in H^2(X,\Z)$ with $w_2(X)\equiv\bar c$
(mod 2) there exists a unique $\Spin^c(4)$-bundle $\hat{P_c}$ with
$P\simeq\qmod{\hat{P_c}}{S^1}$, and $c_1(\det(\hat{P_c}))=c$.
\end{lm}

We denote by $\Sigma_c^{\pm}$ the induced $\U(2)$-vector bundles, and we put
$\Sigma_c:=\Sigma^+_c\oplus\Sigma^-_c$.
\begin{lm} {\rm [1], [11]} The choice of a $\Spin^c(4)$-lift $\hat{P_c}$ of $P$
induces an isomorphism
$$\gamma_+:\Lambda^1\otimes\C\map\Hom_{\C}(\Sigma^+_c,\Sigma^-_c)$$
satisfying the identity
$\gamma_+(u)^*\gamma_+(v)+\gamma_+(v)^*\gamma_+(u)=2g(u,v)\id_{\Sigma^+_c}$
for {\sl
real} cotangent vectors $u,\ v\in\Lambda^1.$
\end{lm}
We define the homomorphisms $\gamma:\Lambda^1\map\End_0(\Sigma_c)$,
$\Gamma:\Lambda^2\map\End_0(\Sigma_c)$ by
$$\gamma(u):=\left(\matrix{0&-\gamma_+(u)^*\cr\gamma_+(u)&0\cr}\right)$$
$$\Gamma(u\wedge v):=\frac{1}{2}[\gamma(u),\gamma(v)]\ ,$$
and we denote by the same symbols
also their $\C$-linear extensions
\linebreak$\Lambda^1\otimes\C\map\End_0(\Sigma_c)$,
and $\Lambda^2\otimes\C\map\End_0(\Sigma_c)$. The homomorphism $\gamma$
defines a map
$\Lambda^1\otimes\Sigma_c\map\Sigma_c$\ ,  called the {\sl Clifford
multiplication}.
The map
$\Gamma$ identifies the bundles $\Lambda^2_{\pm}$ with the bundles of {\sl
trace free
skew-Hermitian} endomorphisms of $\Sigma^{\pm}_c$.

Fix a $\Spin^c(4)$-bundle $\hat{P_c}$ with $P\simeq\qmod{\hat{P_c}}{S^1}$,
and let $L_c:=\det(\hat{P_c})$ be the associated $S^1$-vector bundle. $L_c$
is the
unique unitary line bundle with Chern class $c$.
\begin{re} {\rm [11]} The choice of a $S^1$-connection $a$ in $L_c$ is
equivalent to the
choice of a $\Spin^c(4)$-connection $A$ in $\hat{P_c}$ projecting onto the
Levi-Civita connection.
\end{re}
\begin{dt} The composition
$\Dr_a:A^0(\Sigma_c)\textmap{\nabla_A}A^1(\Sigma_c)\stackrel{\gamma}{\map}A^
0(\Sigma_c)$
is called the Dirac operator associated to the connection $a\in{\cal A}(L_c)$.
\end{dt}
{\bf Notation:} Let  ${\cal A}(L_c)$ be the affine space of $S^1$-connections
in $L_c$. For a connection $a\in{\cal A}(L_c)$, we denote by $F_a\in
A^2(\ad(L_c))=
iA^2$ its curvature, and by $F_a^{\pm}\in i A^2_{\pm}$ the  components of $F_a$
with respect to the orthogonal splitting $A^2=A^2_+\oplus A^2_-$.
Every spinor $\Psi\in A^0(\Sigma^+_c)$ has a conjugate $\bar\Psi\in
A^0(\bar\Sigma^+_c)$, and we can interpret $\Psi\otimes\bar\Psi$ as a Hermitian
endomorphism of $\Sigma^+_c$. Let $(\Psi\otimes\bar\Psi)_0\in
A^0(\End_0(\Sigma^+_c))$ denote the trace-free component of  it.
\vspace{4mm}\\

The {\sl monopole equations} for a pair $(a,\Psi)\in
{\cal A}(L_c)\times A^0(\Sigma^+_c)$ are the equations [18]:
$$\left\{\begin{array}{ll} \Dr_a\Psi=&0\\
\Gamma(F_a^+)=&2(\Psi\otimes\bar\Psi)_0
\end{array}\right.\eqno{(SW)}$$

\begin{pr} {\rm (The Weitzenb\"ock formula [11])}. Let $s$ be the scalar
curvature of
$(X,g)$. Fix a
$\Spin^c(4)$-structure on $X$, and choose a $S^1$-connection $a\in{\cal
A}(L_c)$.
Then the following identity holds on $A^0(\Sigma_c)$ :
$$\Dr_a^2=\nabla_A^*\nabla_A+\frac{1}{2}\Gamma(F_a)+\frac{s}{4}\id_{\Sigma_c
}\ .$$
\end{pr}
\begin{co} Let $\Psi\in A^0(\Sigma^+_c)$. Then
$$\parallel\Dr_a\Psi\parallel^2+\frac{1}{2}\parallel\frac{1}{2}\Gamma(F_a^+)-
(\Psi\bar\Psi)_0)\parallel^2=\parallel\nabla_A\Psi\parallel^2+
\frac{1}{8}\parallel F_a^+\parallel^2+\frac{1}{4}\parallel\Psi\parallel^4+
\frac{1}{4}\int\limits_X s|\Psi|^2.$$
\end{co}
\pf By the Weitzenb\"ock formula we have
$$(\Dr^2_a\Psi,\Psi)=(\nabla_A^*\nabla_A\Psi,\Psi)+
 \frac{1}{2}(\Gamma(F_a^+)(\Psi),\Psi)+\frac{s}{4}(\Psi,\Psi)\ ,$$
since $\Gamma(F_a^-)$ vanishes on $\Sigma_c^+$; integration over $X$ yields:
$$\matrix{\parallel\Dr_a\Psi\parallel^2+\frac{1}{2}\parallel\frac{1}{2}
\Gamma(F_a^+)-
(\Psi\bar\Psi)_0)\parallel^2=\int\limits_X
(\Dr^2_a\Psi,\Psi)+\frac{1}{2}\int\limits_X|\frac{1}{2}\Gamma(F_a^+)-
(\Psi\bar\Psi)_0|^2=\cr
=\parallel\nabla_A\Psi\parallel^2+
\frac{1}{2}\int\limits_X(\Gamma(F_a^+),(\Psi\bar\Psi)_0)+\frac{1}{4}\int\lim
its_X
s|\Psi|^2+\cr+\frac{1}{2}\int\limits_X\frac{1}{4}|\Gamma(F_a^+)|^2-
\frac{1}{2}\int\limits_X(\Gamma(F_a^+),(\Psi\bar\Psi)_0)+
\frac{1}{4}\parallel\Psi\parallel^4\ .} $$

\begin{re} {\rm [18]} If $s\geq 0$ on $X$, then the only solutions $(a,\Psi)$
of
$(SW)$ are pairs $(a,0)$ with $F_a^+=0$.
\end{re}
\section{Seiberg-Witten Invariants}

The {\sl gauge group} ${\cal G}:={\cal C}^{\infty}(X,S^1)$ in the
Seiberg-Witten
theory {\sl is abelian} and acts on ${\cal A}(L_c)\times A^0(\Sigma^+_c)$ by
$(a,\Psi)\cdot f:=(a+f^{-1}df,f^{-1}\Psi)$, letting invariant the set of
solutions
of the equations $(SW)$. We denote by ${\cal W}^g_X(c)$ the moduli space of
solutions of the Seiberg-Witten equations, modulo gauge equivalence. A standard
technique provides a natural structure of finite dimensional {\sl real
analytic}
space in ${\cal W}^g_X(c)$ [6], [5], [16]. The {\sl expected dimension} of this
moduli space is
$$w_c=\frac{1}{4}(c^2-2e(X)-3\sigma(X))\ ,$$
where $e(X)$ and $\sigma(X)$ stand for the Euler characteristic and the
signature
of the oriented manifold $X$. A solution $(a,\Psi)$ is reducible (has
nontrivial
stabilizer) if and only if $\Psi=0$, and then the connection $a$ must be
anti-selfdual. We say that the metric $g$ is $c$-{\sl good} if the $g$-harmonic
representative of the de Rham cohomology class $c_{\rm DR}$ is {\sl not}
anti-selfdual.
If $g$ is $c$-good, then ${\cal W}^g_X(c)$ consists only of irreducible orbits.

Using the same technique as in Yang-Mills theory
([6], [5]), one defines a {\sl gauge invariant perturbation} of the
Seiberg-Witten
equations in order to get smooth moduli spaces of the expected dimension.
For a selfdual
form $\mu\in A^2_+$ we denote by ${\cal W}^{g,\mu}_X(c)$ the moduli space
of solutions
of the perturbed Seiberg-Witten equations
$$\left\{\begin{array}{lll} \Dr_a\Psi&=&0\\
\Gamma(F_a^+ +i\mu)&=&2(\Psi\otimes\bar\Psi)_0
\end{array}\right.\eqno{(SW_{\mu})}$$
We refer to [15] for the following
\begin{lm}\hfill{\break}
1. For every $\mu\in A^2_+$, the moduli space ${\cal W}^{g,\mu}_X(c)$ is
compact.\hfill{\break}
2. There is a dense, second category set of perturbations $\mu\in A^2_+$, for
which the
irreducible part
${\cal W}^{g,\mu}_X(c)^*$  of ${\cal W}^{g,\mu}_X(c)$ is smooth and has the
expected
dimension.\hfill{\break}
3. If $g$ is $c$-good, and $\mu$ is small enough in the $L^2$ topology,
then ${\cal
W}^{g,\mu}_X(c)$  consists only of irreducible orbits, i.e. ${\cal
W}^{g,\mu}_X(c)$=${\cal W}^{g,\mu}_X(c)^*$ .\hfill{\break}
4. Let $g_0$ and $g_1$ be $c$-good metrics which  can be connected
by a smooth path of $c$-good metrics, and let $\varepsilon_i>0$ be small
enough such
that ${\cal W}^{g_i,\mu_i}_X(c)={\cal W}^{g_i,\mu_i}_X(c)^*$ for all
perturbations
$\mu_i$ with $\parallel\mu_i\parallel<\varepsilon_i$. Then  any two
moduli spaces
${\cal W}^{g_i,\mu_i}_X(c)\ ,\  \ i=0, 1$ , with
$\parallel\mu_i\parallel<\varepsilon_i$, which are smooth and have the expected
dimension,  are cobordant.
\end{lm}

The first assertion is a simple consequence of the Weitzenb\"ock formula
and of the
Maximum Principle. The other three assertions follow as in Donaldson theory
by the
Sard theorem for smooth Fredholm maps, and by transversality arguments [5].
Note that
in 4. we mean cobordism between {\sl non-oriented} compact smooth
manifolds. A more
delicate analysis of the monopole equations [18] shows  that, in fact, the
moduli
spaces ${\cal W}^{g,\mu}_X(c)^*$  come with
natural orientations, as soon as they are smooth and have the expected
dimension, and
that the conclusion in 4. holds for the oriented moduli spaces.
\vspace{3mm}\\

The Seiberg-Witten theory  provides strong differentiable
invariants using only moduli spaces of dimension 0. Let $c$ be an integral
lift of
$w_2(X)$, with $w_c=0$, i.e.
$$c^2=2e(X)+3\sigma(X)\ .$$

Such a lift is called an {\sl almost canonical class}, since the condition
\hbox{$w_c=0$} is equivalent to the existence of an almost complex structure on
$X$ with first Chern class $c$ [10], [16].

Now fix an almost canonical class $c$, choose a $c$-good metric $g$, and a
small,
sufficiently general perturbation $\mu$. Then ${\cal W}^{g,\mu}_X(c)={\cal
W}^{g,\mu}_X(c)^*$ is compact, smooth of the expected dimension 0, and its
bordism
class is independent of $\mu$. Let
$n_c^g:=|{\cal W}^{g,\mu}_X(c)|$ mod 2 be the number of points modulo 2 of this
moduli space. Lemma 2.1 implies that
$n_c^g$ is also independent of $g$ if any two $c$-good metrics can be
connected by a
smooth 1-parameter family of $c$-good metrics.

 The numbers  $n_c:=n_c^g$ associated to such almost canonical
classes are called the   mod 2-Seiberg-Witten  invariants, and the classes
$c$ with
$n_c\ne 0$ are then called mod 2-Seiberg-Witten classes of index 0. By
definition
they are  differentiable invariants, in the following sense: If $f:X'\map
X$ is an
orientation-preserving diffeomorphism, and for an almost canonical class $c$ of
$X$ the Seiberg-Witten invariant $n_c$ is well defined, then $f^*(c)$ has the
same property, and $n_{f^*(c)}=n_{c}$.
\begin{re}
Let $c$ be an almost canonical class of $X$. \hfill{\break}
1. If $c^2\geq 0$ and $c_{\rm DR}\ne 0$ , then {\sl any} Riemannian metric
on $X$ is
$c$-good.\hfill{\break} 2. If $b_2^+\geq 2$, then any two $c$-good metrics
can be
connected  by a smooth path of $c$-good metrics.
\end{re}

Therefore, if one of the two conditions above is satisfied, then the
\linebreak\hbox{mod
2-Seiberg-Witten} invariant $n_c$ is well-defined.
\vspace{3mm}

In the case $b_2^+=1$, invariants can still be defined, but the dependence of
$n_c^g$ on the metric $g$ must be taken into account: In the real vector space
$H^2_{\rm DR}(X)$ , consider the positive cone
$${\cal K}=\{u\in H^2_{\rm DR}(X) \ | u^2>0\}\ .$$
Fix a {\sl non-vanishing} cohomology class $k\in H^2_{\rm DR}(X)$ with
$k^2\geq 0$. The
cone
${\cal K}$ splits as the disjoint union of its connected components ${\cal
K}_{\pm}$,
where
$${\cal K}_{\pm}:=\{u\in{\cal K} \ |\ \pm u\cdot k>0\}\ .$$

If $c$ is
an almost canonical class, let $c^{\bot}$ be the hyperplane
$$c^{\bot}:=\{u\in H^2_{\rm DR}(X)\ |\ c\cdot u=0\}$$

If $c^{\bot}$ meets ${\cal K}_+$, then the intersection $c^{\bot}\cap{\cal
K}_+$ is
called the  {\sl  wall} of type $c$, and the two components of ${\cal
K}_+\setminus
c^{\bot}$ are called {\sl chambers} of type $c$.
%
%
For every Riemannian metric $g$ on $X$, let $\omega_g$ be a generator of
the real line
of $g$-harmonic selfdual 2-forms, such that $[\omega_g]\in{\cal K}_+$ .
Then the ray
$\R_{>0}[\omega_g]\subset {\cal K}_+$ depends smoothly on the metric $g$.
The property of a
metric to be $c$-good has the following simple geometric interpretation:
\begin{re}\hfill{\break} Suppose $b_2^+(X)=1$. Then:\hfill{\break}
1. The metric $g$ is $c$-good iff the ray $\R_{>0}[\omega_g]$ does not lie in
the wall $c^{\bot}\cap{\cal K}_+$. \hfill{\break}
2. If $g_0$, and $g_1$ are $c$-good metrics, then $n_c^{g_0}=n_c^{g_1}$ iff
the two
rays  $\R_{>0}[\omega_{g_i}]$  belong to the same chamber of type $c$.
\end{re}

The first assertion follows immediately from the definition. The second needs a
careful analysis of a 1-parameter family of 0-dimensional smooth
moduli spaces  ${\cal W}^{g_t}_X(c)$ around the value of the parameter $t$ for
which the ray
$\R_{>0}[\omega_{g_t}]$ crosses the wall $c^{\bot}\cap{\cal K}_+$ (see
[18], [15]).

\section{Monopoles on K\"ahler surfaces}

Let $(X,J,g)$ be an almost complex 4-manifold endowed with a Hermitian
metric $g$. The
almost complex structure $J$ defines a reduction of the structure group of the
tangent bundle $T_X$ of $X$ from $\SO(4)$ to $\U(2)$. In particular, we get
a {\sl
canonical}
$\Spin^c(4)$-structure on \ $X$ \ via \ the \ canonical \ lifting \
\ $l:\U(2)\map\Spin^c(4)$ [11]. Let $\omega_g$ be the K\"ahler form of $g$.
\begin{lm} {\rm [11]} The canonical $\Spin^c$-structure of an almost
complex Hermitian
4-manifold has the following properties:\hfill{\break}
1. There are canonical identifications
$\Sigma^+=\Lambda^{00}\oplus\Lambda^{02}$,
$\Sigma^-=\Lambda^{01}$.\hfill{\break}
2. Via these identifications, the map
$\Gamma:\Lambda^2_+\otimes\C\map\End_0(\Sigma^+)$ is given by:
$$\Lambda^{20}\oplus\Lambda^{02}\oplus\Lambda^{00}\omega_g\ni(\lambda^{20},\
lambda^{02},
f\omega_g)\textmap{\Gamma}2\left[\matrix{-if&-*(\lambda^{20}\wedge\cdot)\cr
\lambda^{02}\wedge\cdot&if\cr}\right]\in\End_0(\Lambda^{00}\oplus\Lambda^{02
})\ .$$
\end{lm}

Suppose now that $(X,J,g)$ is a K\"ahler surface. This means that $J$ is
integrable,
and $\omega_g$ is closed (or equivalently, $J$ is Levi-Civita parallel). In
particular
the holonomy group of the Levi-Civita connection also reduces to $\U(2)$,
and the
splittings $\Lambda^p\otimes\C=\bigoplus\limits_{i+j=p}\Lambda^{ij}$ are
Levi-Civita
parallel. We get a $\U(2)$-connection in the holomorphic tangent bundle ${\cal
T}_X=T^{10}_X\simeq\Lambda^{01}$, which coincides with the {\sl Chern
connection} of
this bundle, i.e. with the unique connection compatible with the
holomorphic structure
and the Hermitian metric. The induced connection $c_0$ in the line bundle
$K_X^{\vee}=\det({\cal T}_X)\simeq\Lambda^{02}$ also coincides  with the Chern
connection of this Hermitian holomorphic line bundle.

Every other $\Spin^c(4)$ structure $\hat{P_c}\map P$ on $(X,g)$ has as
spinor bundle
$$\Sigma_c=\Sigma\otimes M\ ,$$
where $M$ is a differentiable $S^1$-bundle with $2c_1(M)+c_1(K_X^{\vee})=c$.
(For a simply connected manifold $X$, $M$ is well defined up to
isomorphy  by this condition.) $S^1$-connections in
$\det(\Sigma^{\pm}_c)=K_X^{\vee}\otimes M^{\otimes 2}$ correspond to
$S^1$-connections
in $M$. Given $b\in{\cal A}(M)$, the curvature of the corresponding connection
$a\in{\cal A}(K_X^{\vee}\otimes M^{\otimes 2})$ is $F_a=F_{c_0}+2F_b$.

A half-spinor $\Psi\in A^0(\Sigma^+\otimes M)$ can be written as
$$\Psi=\varphi+\alpha\ ,\ \ \varphi\in A^0(M)\ ,\ \ \alpha\in
A^{02}(M)\ .$$

We put $J(M):=c_1(\Sigma^+\otimes M)\cup[\omega_g]$.
\begin{pr} Let $(X,g)$ be a K\"ahler surface with Chern connection $c_0$ in
$K_X^{\vee}$, $M$ a differentiable $S^1$-bundle with $J(M)<0$. A pair
$(b,\varphi+\alpha)\in {\cal A}(M)\times\left(A^0(M)\oplus A^{02}(M)\right)$
solves the
monopole equations iff:
$$
\begin{array}{l}F_b^{20}=F_b^{02}=0\  \\
\alpha=0\ ,\ \ \bar\partial_b(\varphi)=0 \ \\
i\Lambda F_b+\frac{1}{2}\varphi\bar\varphi+\frac{s}{2}=0\ .\end{array}\
\eqno{(*)}$$
\end{pr}
\pf The pair $(b,\varphi+\alpha)$ solves the equations $(SW)$ iff the
corresponding pair $(a,\varphi+\alpha)$ satisfies
$$\begin{array}{ll}F_a^{20}&=-\varphi\otimes\bar\alpha\\
F_a^{02}&=\ \alpha\otimes\bar\varphi\\
\bar\partial_b(\varphi)&=\ i\Lambda\partial_b(\alpha) \\
i\Lambda
F_a&=-\left(\varphi\bar\varphi-*(\alpha\wedge\bar\alpha)\right).\end{array}\
$$
By Corollary 1.7 it follows that  $(b,\varphi+\alpha)$ solves $(SW)$ iff
$(b,\varphi-\alpha)$ does (Witten's trick).
Therefore
$\varphi\otimes\bar\alpha=\alpha\otimes\bar\varphi=0$, hence
$F_a^{20}=F_a^{02}=0$, and
$\varphi$ or $\alpha$ must vanish. Integrating the equation $i\Lambda
F_a=-\left(\varphi\bar\varphi-
*(\alpha\wedge\bar\alpha)\right)$ over $X$, we find:
$$J(M)=(2c_1(M)-c_1(K_X))\cup[\omega_g]=\int\limits_X\frac{i}{2\pi}
F_a\wedge\omega_g=
\frac{1}{8\pi}\int\limits_X(-|\varphi|^2+|\alpha|^2)\ ,$$
hence  $\alpha=0$ if $J(M)<0$.
\qed

The above proposition must be interpreted as follows: If $J(M)<0$, then the
solutions
of the monopole equations $(SW)$ are the pairs $(b,\varphi)\in{\cal A}(M)\times
A^0(M)$, such that
$b$ is the Chern connection of a holomorphic structure in $M$, $\varphi$ is a
holomorphic section, and the mean curvature $i\Lambda F_b$ of $b$ satisfies
the {\sl
generalized vortex equation} [16], [3], [4], [9]
$$i\Lambda F_b+\frac{1}{2}\varphi\bar\varphi+\frac{s}{2}=0\ .\eqno{(V_s)}$$

Moreover, every {\sl infinitesimal deformation} of a solution of the form
$(b,\varphi)$, $\varphi\ne 0$ of the monopole equation still vanishes in the
$\alpha$-direction. Therefore ${\cal W}_X^g(c)$ can be identified (as {\sl real
analytic space}) with the moduli space of pairs
$(b,\varphi)$ satisfying the above conditions, modulo the gauge group ${\cal
C}^{\infty}(X,S^1)$ of unitary automorphisms of $M$. Under the assumption
$J(M)<0$,
the action of the gauge group is free on the space of solutions, because
any solution
$(b,\varphi)$ has a non-vanishing section $\varphi$.

Alternatively, let ${\cal M}$ be a holomorphic line bundle with differentiable
support $M$, and $\varphi$ a holomorphic section of ${\cal M}$. For a
Hermitian metric
$h$ in ${\cal M}$, we denote by $F_h$ the curvature of the associated Chern
connection,
and we consider the following equation for $h$:
$$i\Lambda F_h+\frac{1}{2}\varphi\bar\varphi^h+\frac{s}{2}=0\ . \eqno{(V'_s)}$$

Standard arguments (see for instance [16], [9]) show that the problem of
classifying
the solutions $(b,\varphi)$ of $(*)$  modulo {\sl unitary automorphisms of
$M$} is
equivalent to the problem of classifying  those pairs $({\cal M},\varphi)$
modulo
{\sl holomorphic isomorphisms}, for which the equation $(V'_s)$ has a solution.
\begin{pr}
Let $(X,g)$ be a compact K\"ahler surface, $({\cal M},\varphi)$ a
holomorphic line
bundle with a {\sl non-vanishing} holomorphic section $\varphi\in
H^0(X,{\cal M})$.
${\cal M}$ admits a metric
$h$ satisfying the equation $(V'_s)$ iff
$$c_1({\cal M})\cup[\omega_g]<\frac{1}{2}c_1(K_X)\cup[\omega_g]\ .$$
\end{pr}
\pf (cf. [3]) Fix a background metric $h_0$; any other metric $h$ has the form
$h=e^{2u}h_0$, with $u\in A^0$ a smooth function.The vortex equation $(V'_s)$
translates into
$$\Delta u+\frac{1}{2}|\varphi|^2_{h_0} e^{2u}+(i\Lambda
F_{h_0}+\frac{s}{2})=0\ .\eqno{(1)}$$
Set $q:=\int\limits_X(i\Lambda F_{h_0}+\frac{s}{2})=2\pi(c_1({\cal
M})-\frac{1}{2}c_1(K_X)\cup[\omega_g]$, and choose
$v\in A^0$ with
$$-\Delta v=(i\Lambda F_{h_0}+\frac{s}{2})-q \ .$$
Define $w:=2(u-v)$. Then $(1)$ is equivalent to the following equation in $w$:
$$\Delta w+(|\varphi|^2_{h_0} e^{2v})e^w +2q=0\ .\eqno{(2)}$$

Integrating over $X$, we see that if (2) has solutions, then
$q$ must be negative. On the other hand, by a well known result of Kazdan
and Warner
[3],  (2) has a unique solution if $q<0$.
\qed
\begin{th} {\rm [18], [16]}
Let $(X,g)$ be a simply connected K\"ahler surface,   $c\in H^2(X,\Z)$ with
$c\equiv c_1(K_X)$ mod 2, and  $\pm c\cup[\omega_g]<0$.
\hfill{\break}
1. \ If $c\ \not\in\ \NS(X)$, then ${\cal W}_X^g(c)=\emptyset$ .\hfill{\break}
2. Suppose $c\in\NS(X)$. Then there is a natural real analytic isomorphism
\linebreak ${\cal W}_X^g(c)\simeq\P(H^0(X,{\cal M}))$, where ${\cal M}$ is
the (unique,
up to isomorphy) holomorphic line bundle with
$c_1(K_X^{\vee}\otimes{\cal M}^{\otimes 2})=\pm c$. \hfill{\break}
3. ${\cal W}^g_X(c)$ is always smooth.\ Let $D$ be the divisor of a
nontrivial section in
${\cal M}$. Then ${\cal W}^g_X(c)$ has the expected dimension iff
\hbox{$h^1({\cal O}_X(D)|_D)=0$}.
\end{th}

\section{Rationality of complex surfaces}

A compact complex surface is {\sl rational} iff its field of meromorphic
functions is
isomorphic to $\C(u,v)$. Such a surface is always simply connected and has
$b_2^+=1$
[2]. The following result has been
has been announced by R. Friedman and Z. Qin [8]. Whereas their proof uses
Donaldson theory and vector bundles techniques, our proof uses the new
Seiberg-Witten invariants, and our interpretation of these invariants in
terms of linear systems.
\begin{th} {\rm [17]}
A complex surface $X$ which is diffeomorphic to a rational surface is rational.
\end{th}
\pf The proof consists of the following three steps:

1.Any rational surface $X_0$ admits a
{\sl Hitchin metric} [12], i.e. a K\"ahler metric
$g_0$ with positive {\sl total scalar curvature}. This condition can be
written as
$c_1(K_{X_0})\cup[\omega_{g_0}]<0$.

Let $c$ be any integral lift of $w_2(X_0)$, such that $g_0$ is $c$-good,
i.e. such that
the moduli space ${\cal W}^{g_0}_{X_0}(c)$ contains no reducible solutions.
Since $p_g(X_0)=0$, $c$ has always type (1,1), and $g_0$ is $c$-good iff
$c\cup[{\omega_{g_0}}]\ne 0$.

We assert
that ${\cal W}^{g_0}_{X_0}(c)$ is then empty, and in particular, all
Seiberg-Witten
invariants $n_c^{g_0}$ computed with respect to this metric vanish.

Indeed, let ${\cal M}$ be the holomorphic line bundle defined in Theorem
3.4. If the
moduli space $\P(H^0(X_0, {\cal M}))$ was not empty, then
$$c_1({\cal M})\cup[\omega_{g_0}]\geq 0 \ .\eqno{(1)}$$
But we have
$$0>\pm c\cup[\omega_{g_0}]=(2c_1({\cal M})-c_1(K_{X_0}))\cup[\omega_{g_0}]\
,$$
hence, by (1)
$$0\leq 2c_1({\cal M})\cup[\omega_{g_0}]<c_1(K_{X_0})\cup[\omega_{g_0}]\ ,$$
which contradicts the assumption on the total scalar curvature of $g_0$.
\vspace{2mm}\\

2. Let now $X$ be a simply connected projective surface with ${\rm kod}(X)>
0$. We may
suppose that $X$ is the blow up in $k$ {\sl distinct} points of its minimal
model
$X_{\min}$.  Denote by
$\sigma:X\map X_{\min}$ the contraction to the minimal model, and by
$E=\sum\limits_{i=1}^k E_i$ the exceptional divisor. Fix an ample divisor
$H_{\min}$ on
$X_{\min}$,  set
$H_n:=\sigma^*(n H_{\min})-E$, and for $n\gg 0$ choose a K\"ahler metric
$g_n$ on
$X$ with $[\omega_{g_n}]=c_1(H_n)$. Given $I\subset\{1,\dots,k\}$, define
$$\begin{array}{l}E_I:=\sum\limits_{i\in I} E_i\\ c_I:=2c_1(E_I)-c_1(K_X)\\
\bar I:=\{1,\dots,k\}\setminus I \ .\end{array}$$
Since $c_I$ is an almost canonical class, the expected dimension of the
corresponding Seiberg-Witten moduli space is 0. For
$n\gg 0$ we get
$c_I\cup [\omega_{g_n}]<0$, and Theorem 3.4 gives
$${\cal W}^{g_n}_X(c_I)\simeq\{E_I\}\ .$$
Therefore ${\cal W}^{g_n}_X(c_I)$ consists of a single smooth point, and
$$n_{c_I}^{g_n}= 1\ {\rm mod}\ 2\ .\eqno{(2)}$$

3. Suppose now that there is an orientation-preserving diffeomorphism
$f:X\map X_0$, where
$X$ is projective surface with ${\rm kod} X\geq 0$. Since $X$ must have
$p_g(X)=0$, and
$\pi_1(X)=\{1\}$, it follows that, in fact, ${\rm kod} X>0$. Let
$g=f^*(g_0)$ denote the pull-back of a Hitchin metric to
$X$; clearly
$$n_{c_I}^g=0 \eqno{(3)}$$
for all $I\subset\{1,\dots,k\}$ such that $g$ is $c_I$-good.

We will now derive a contradiction in the following way: Using the
Enriques-Kodaira
classification of surfaces, it easy to see that the de Rham cohomology class
$ k_{\min}:=\sigma^*(c_{1,{\rm DR}}(K_{\min}))$ is non-trivial and satisfies
the condition
$k^2_{\min}\geq 0$. Therefore we can consider the upper positive cone
$${\cal K}_+:=\{u\in H^2_{\rm DR}(X)\ |\ u^2>0, \ u\cdot k_{\min}>0\}\ .$$

Clearly $[\omega_{g_n}]$ belongs to ${\cal K}_+$. We choose a harmonic
$g$-selfdual form
$\omega_g$, with
$[\omega_g]\in{\cal K}_+$.
\vspace{3mm}\\
{\bf Claim:} {\sl The rays $\R_{>0}[\omega_g]$ and $\R_{>0}[\omega_{g_n}]$
belong
either to the same chamber of type $c_I$ or to the same chamber of type
$c_{\bar I}$.}
\vspace{1mm}\\
\pf If not, then, since $c_I\cup [\omega_{g_n}]<0$,  we get
$[\omega_g]\cdot c_I\geq0$
and
$[\omega_g]\cdot c_{\bar I}\geq 0$. Write
$$[\omega_g]=\sum\limits_{i=1}^k\lambda_i E_i+\sigma^*[\omega]\ ,$$
with $[\omega]\in
H^2_{\rm DR}(X_{\min})$. Then
$$\begin{array}{l} -\sum\limits_{i\in I}\lambda_i+\sum\limits_{j\in\bar
I}\lambda_j-[\omega]\cdot [K_{\min}]\geq 0\\
-\sum\limits_{j\in\bar I}\lambda_j+\sum\limits_{i\in
I}\lambda_i-[\omega]\cdot [K_{\min}]\geq 0\ .\end{array}$$
Adding these inequalities we find $[\omega]\cdot[K_{\min}]\leq 0$. But
$[\omega]\cdot [K_{\min}]=[\omega_g]\cdot k_{\min}>0$, because
$[\omega_g]\in{\cal K}_+$. This contradiction proves the claim.

\vspace{3mm}

It follows that either $g$ and $g_n$ are both $c_I$-good and
$n_{c_I}^g=n_{c_I}^{g_n}$, or  $g$ and $g_n$ are both $c_{\bar I}$-good and
$n_{c_{\bar I}}^g=n_{c_{\bar I}}^{g_n}$. This gives now a contradiction
with (2) and
(3).
\qed

Together with the results of Friedman and Morgan [7], we have:
\begin{th}{\rm (The Van de Ven conjecture [19])} The Kodaira dimension of
complex
surfaces is a
${\cal C}^{\infty}$-invariant.
\end{th}
%
%
%
%
{\bf Remark:} It is possible to couple the Seiberg-Witten equations to
connections
in unitary bundles. The solutions of these coupled Seiberg-Witten equations
over
K\"ahler surfaces again have a purely complex-geometric interpretation
[16]: The moduli
space of solutions can be identified---via generalized vortex equations--- with
moduli spaces of stable pairs [13], [4]. This construction could lead to
new invariants
which might be nontrivial for K\"ahler surfaces wit $p_g=0$.
\vspace{0.8cm}\\
\parindent0cm
\centerline {\Large {\bf Bibliography}}

\vspace{0.5cm}
1. Atiyah M., Hitchin N. J., Singer I. M.: {\it Selfduality in
four-dimensional Riemannian geometry}, Proc. R. Lond. A. 362, 425-461 (1978)

2. Barth, W., Peters, C., Van de Ven, A.: {\it Compact complex surfaces},
Springer Verlag (1984)

3. Bradlow, S. B.: {\it Vortices in holomorphic line bundles over closed
K\"ahler manifolds}, Comm. Math. Phys. 135, 1-17 (1990)

4. Bradlow, S. B.: {\it Special metrics and stability for holomorphic
bundles with global sections}, J. Diff. Geom. 33, 169-214 (1991)

5. Donaldson, S.; Kronheimer, P. B.: {\it The Geometry of four-manifolds},
Oxford Science Publications (1990)

6. Freed, D. S.; Uhlenbeck,  K. K.: {\it Instantons and Four-Manifolds},
Springer Verlag (1984)

7. Friedman, R., Morgan, J.W.: {\it Smooth 4-manifolds and Complex Surfaces},
Springer Verlag  3. Folge, Band 27 (1994)

8.  Friedman, R., Qin, Z.: {\it On complex surfaces diffeomorphic to
rational surfaces}, Preprint (1994)

9. Garcia-Prada, O.: {\it Dimensional reduction of stable bundles, vortices
and stable pairs}, Int. J. of Math. Vol. 5, No 1, 1-52 (1994)

10. Hirzebruch, F., Hopf H.: {\it Felder von Fl\"achenelementen in
4-dimensionalen 4-Mannigfaltigkeiten}, Math. Ann. 136, (1958)

11. Hitchin, N.: {\it  Harmonic spinors}, Adv. in Math. 14, 1-55 (1974)

12. Hitchin, N.: {\it  On the curvature of rational surfaces}, Proc. of Symp.
in Pure Math., Stanford, Vol. 27 (1975)

13. Huybrechts, D.; Lehn, M.: {\it Stable pairs on curves and surfaces},
 J. Alg. Geometry, (1995)

14. Kobayashi, S.: {\it Differential geometry of complex vector bundles},
Princeton University Press, (1987)

15. Kronheimer, P., Mrowka, T.: {\it The genus of embedded surfaces in the
projective plane}, Preprint (1994)

16. Okonek, Ch.; Teleman A.: {\it The Coupled Seiberg-Witten Equations,
Vortices, and Moduli Spaces of Stable Pairs}, Preprint, January, 13-th 1995

17. Okonek, Ch.; Teleman A.: {\it Seiberg-Witten invariants and the Van de
Ven conjecture}, Preprint, February, 8-th 1995

18. Witten, E.: {\it Monopoles and four-manifolds}, Mathematical Research
Letters 1, 769-796 (1994)

19. Van de Ven, A,: {\it On the differentiable structure of certain algebraic
surfaces}, S\'em. Bourbaki ${\rm n}^o$ 667, Juin (1986)
\vspace{0.5cm}\\
Authors addresses:\\
Mathematisches Institut, Universit\"at Z\"urich,\\
Winterthurerstrasse 190, CH-8057 Z\"urich\\
e-mail:okonek@math.unizh.ch

\ \ \ \ \ \  \ \ \ teleman@math.unizh.ch

\end{document}